\documentclass[10pt]{amsart}

\usepackage{amsmath}

\usepackage{amssymb}

\usepackage{graphicx}

\newtheorem{thm}{Theorem}[section]
\newtheorem{prop}[thm]{Proposition}
\newtheorem{lemma}[thm]{Lemma}

\theoremstyle{definition}
\newtheorem{defn}[thm]{Definition}

\theoremstyle{remark}
\newtheorem{remark}[thm]{Remark}

\numberwithin{equation}{section}

\def\R{\mathbb{R}}

\def\Z{\mathbb{Z}}

\def\P{\mathbb{P}}

\def\F{\mathcal{F}}

\def\hom{{\mathrm{Hom}}}
\def\g{\mathfrak{g}}
\def\a{\mathfrak{a}}
\def\b{\mathfrak{b}}
\def\m{\mathfrak{m}}
\def\n{\mathfrak{n}}

\def\sl{\mathfrak{sl}}

\def\SO{\mathrm{SO}}
\def\SL{\mathrm{SL}}

\def\I{\mathcal{I}}
\def\A{\mathcal{A}}

\begin{document}

\title[Tableaux over Lie algebras]{Tableaux over Lie algebras,
integrable systems, and classical surface theory}

\author{Emilio Musso} 
\address{(E. Musso) Dipartimento di Matematica Pura ed Applicata,
Universit\`a degli Studi dell'Aquila, Via Vetoio, I-67010
Coppito (L'Aquila), Italy} \email{musso@univaq.it}

\author{Lorenzo Nicolodi}
\address{(L. Nicolodi) Di\-par\-ti\-men\-to di Ma\-te\-ma\-ti\-ca,
Uni\-ver\-si\-t\`a degli Studi di Parma, Viale G. P. Usberti 53/A,
I-43100 Parma, Italy} \email{lorenzo.nicolodi@unipr.it}

\thanks{Authors partially supported by the MIUR project
\textit{Propriet\`a Geometriche delle Variet\`a Reali e
Complesse},
%by the GNSAGA of the INdAM,
and by The European Contract
Human Potential Programme, Research Training Network
HPRN-CT-2000-00101 (EDGE)}


%%%%%%%%%%%
%% 11/11/2004 modifiche nella dimostrazione della involutivita' 
%% del sistema associato al tableau; modificato ulteriormente il 22/11/04; 
%% last revised 23/12/04; revised according to the referee reports 17/06/05 
%% and again 30/06/05; last revised 06/07/05; accepted 08/11/05.
%%%%%%%%%%%

%\dedicatory{}


\begin{abstract}
Starting from suitable tableaux over finite dimensional Lie algebras,
we provide a scheme for producing involutive linear Pfaffian systems 
related to various classes of submanifolds in homogeneous spaces which constitute 
integrable systems.
These include isothermic surfaces, Willmore surfaces, and other classical 
soliton surfaces.
Completely integrable equations such as the
$G/G_0$-system of Terng and the curved flat system of Ferus--Pedit
may be obtained as special cases of this construction.
Some classes of surfaces in projective differential geometry whose Gauss--Codazzi
equations are associated with tableaux over $\mathbf{sl}(4,\mathrm{R})$ are discussed.
\end{abstract}

\maketitle

\section{Introduction}\label{s:intro}

\noindent It is well known that various completely integrable nonlinear PDEs (soliton equations)
arise in differential geometry as compatibility conditions for the linear equations satisfied by 
frames adapted to submanifolds in higher dimensional homogeneous spaces. 
In this geometric setting, the adapted frames and the local differential invariants of a submanifold
may be described as integral manifolds of an appropriate exterior differential system (EDS).
The theory of EDSs has been proven to be fruitful in the study of certain
soliton equations, especially as concerns their B\"acklund transformations and conservation laws.
In this regard, we only mention the early work of Estabrook and Wahlquist
on the Korteweg--de Vries equation and the Schr\"odinger equation with a cubic nonlinearity  
\cite{EW1, EW2}, the work of Sasaki, Chern and his collaborators on special classes
of nonlinear evolution equations \cite{Sa, Ce1, Ce2}, and more recently 
that of Bryant, Terng and Wang on the curved flat equations \cite{Br, TW03}. 

\vskip0.2cm

%%%%%%%
The purpose of this paper is to present a unifying approach to a number of 
different classes of known integrable surfaces and submanifolds in homogeneous 
spaces.
These all ultimately derive from a special class of involutive Pfaffian
differential systems.
As such, the EDS point of view offers a way of better understanding the geometry
of these submanifolds and yet another perspective in the search for new classes
of integrable geometries.

\vskip0.2cm

We provide a scheme for producing linear Pfaffian systems in 
involution\footnote{In the
sense of Cartan-K\"ahler theory \cite{BCGGG}.}
starting from suitable algebraic tableaux over finite dimensional Lie algebras.
%%%%%%%%%%%%%%
The interest in this construction is that several classes of submanifolds
in homogeneous spaces which constitute integrable systems are related to Pfaffian systems 
which fit into this scheme.
The solutions of these systems contain, as special cases,
the solutions to well-known completely integrable equations
such as the $n$-dimensional system associated with a symmetric space $G/G_0$ of rank $n$, 
the so-called $G/G_0$-system defined by Terng \cite{Ter97}, and the (gauge equivalent) 
system for curved flats in symmetric spaces introduced by Ferus and Pedit \cite{FP1}.
Geometrically, besides the examples of isometric immersions of space forms into space forms,
of isothermic surfaces and of $L$-isothermic surfaces,
all related to curved flat systems \cite{FP2, BHPP, BDPT02, MN2, MN4}, other important examples 
include Willmore surfaces and their counterparts in Laguerre geometry
\cite{M1, MN96},  
Godeaux--Rozet, Demoulin and asymptotically-isothermic surfaces in projective differential 
geometry, as well as their analogues in Lie sphere geometry \cite{BH, Fe1, Fe3, FS}.

\vskip0.2cm

The common structure of the above examples is the concept of tableau over a Lie 
algebra. It takes up the familiar notion of an involutive tableau
in the theory of EDSs (cf. \cite{Ma54, SS, BCGGG}) and adapts it to the purpose.
Although its simple formulation, we are not aware of it having appeared earlier 
in the literature.
%%%%%%%%%%%%%%%
We illustrate the general procedure to associate with each tableau over a Lie 
algebra a Pfaffian differential system, and show that the property of being involutive 
for such systems can be directly detected from the properties of the corresponding tableaux. 
In the semisimple case, this
construction applies to the $G/G_0$-system and the curved flat system.
At the same time, we start the program of identifying the submanifold geometry associated with 
these systems/tableaux. This amounts to finding submanifolds in some homogeneous space
whose integrability conditions are given by the Pfaffian differential system associated with 
a given tableau. 
We carry out this program for surfaces in projective differential geometry 
whose Gauss--Codazzi equations correspond to various tableaux over $\sl(4,\R)$. 

\vskip0.2cm

The paper is organized as follows. In Section \ref{s:tableaux}, we introduce the notion of a 
linear tableau over a finite dimensional Lie algebra and examine its basic properties, 
discussing the special case of semisimple Lie algebras.
In Section \ref{s:system-tableaux}, starting from a tableau $\A$ over a real Lie algebra 
$\mathfrak{g}$ with Lie group $G$, we construct
a canonical Pfaffian exterior differential system with 
independence condition $(\mathcal{I},\omega)$ on $Y = G\times \A$, and prove that
this system is involutive in the sense of Cartan--K\"ahler theory. In particular, the Cartan
characters of a regular flag are computed in terms of the characters of the corresponding 
tableau. Next, we illustrate the construction of this Pfaffian exterior differential system
in the case of the $G/G_0$-system of Terng.
In Section \ref{s:projective}, we give a brief review of the moving frame method for surfaces 
in real projective space and explain the relation between the exterior differential systems
defined by tableaux over $\sl(4,\R)$ and
the equations for various classes of surfaces in projective differential geometry.
We discuss the cases of Godeaux--Rozet, Demoulin, and asymptotically-isothermic surfaces.
The corresponding tableaux are studied and their characters are computed.
%%%%%%%
\vskip0.2cm

The next step would be to develop a systematic study of tableaux over Lie
algebras. This will help to answer a number of questions which 
naturally arise.
For instance, one can ask which tableaux generate integrable differential systems
and whether it is possible to formulate additional restrictions on the tableaux
to guarantee the integrability.
These and other questions will be the object of future investigation.
We conclude by observing that
the above setting is also the natural one for studying invariant variational problems, 
their Euler--Lagrange systems and conservation laws, in the spirit
of Bryant, Griffiths and Grossman (cf. \cite{BGG} and the notion of 
multi-contact manifold therein).

\vskip0.2cm

The basic reference on exterior differential systems has been \cite{BCGGG}
and our notation is consistent with this reference.
The summation convention over all pairs of repeated indices in a product will be used 
throughout except when explicitly stated otherwise.

\vskip0.2cm
 We thank the referees for useful comments and suggestions.

\section{Tableaux over Lie algebras}\label{s:tableaux}

\subsection{Basic definitions}

Let $(\g, [\,,])$ be a finite dimensional Lie algebra (over any field of characteristic zero). If
$\a, \b$ are vector subspaces of $\g$ such that $\a\cap\b= \{0\}$, let $\A$ be a linear subspace of
$\hom(\a,\b)$ and denote by
\[
 \rho : \hom(\a,\A) \to \b\otimes \Lambda^2(\a^\ast)
  \]
the linear map defined by
\begin{equation}\label{ro}
 \rho(F) (A_1,A_2) = \frac{1}{2}\left(F(A_1)(A_2) - F(A_2)(A_1)\right),
  \end{equation}
for all $F\in \hom(\a,\A)$ and $A_1,A_2 \in \a$. The {\it first prolongation} $\A^{(1)} \subset \hom(\a,\A)$
of the subspace $\A$ is the kernel of $\rho$. If
$\hom(\a,\b)$ is identified with $\b\otimes \a^\ast$, then $\A^{(1)} =
\left(\A\otimes \a^\ast\right)
\cap \left(\b\otimes S^2(\a^\ast)\right)$, where $S^2(\a^\ast)$ is the 2nd symmetric power of $\a^\ast$.

\vskip0.3cm

For any flag $\ell\in \F(\a))$, $\ell: (0)\subset \a_1 \subset\cdots\subset \a_k =\a$, we set
\[
\A_j(\ell) = \{Q\in \A\, : \, Q_{|\a_j} =0\} \quad (j=1,\dots,k).
  \]
The subspaces $\A_j$ give a filtration
\begin{equation}\label{filtrazione}
 (0) = \A_k(\ell)\subseteq \A_{k-1}(\ell) \subseteq\cdots\subseteq \A_0(\ell) = \A.
  \end{equation}
A flag $\ell$ is {\it generic} with respect to $\A$ if the $\text{dim}\,\A_j(\ell)$ 
are a minimum, i.e.,
\[
 \text{dim}\,\A_j(\ell) = \min\{\text{dim}\,\A_j(\tilde{\ell})\, : \, \tilde{\ell} \in \F\}.
  \]
The set of generic flags is an open and dense subset of the flag manifold $\F(\a)$.

\vskip0.3cm

The {\it characters} of $\A$ are the non-negative integers $s'_j$ defined inductively by
\[
 s'_1 +\cdots+s'_j= \text{dim}\, \A - \text{dim}\,\A_j(\ell) \quad (j=1,\dots,k),
  \] where $\ell$ is a generic
flag with respect to $\A$. With this formulation, one can establish the following inequality 
(cf. \cite{Ma54, SS, BCGGG})
\begin{equation}\label{inequality}
 \text{dim}\, \A^{(1)} \leq s'_1 +2s'_2 + \cdots + ks'_k,
  \end{equation}
where
\[
 s'_1 +2s'_2 + \cdots + ks'_k = \text{dim}\, \A +\text{dim}\, \A_1 +\cdots + \text{dim}\, \A_{k-1},
  \]
If equality holds in \eqref{inequality}, then the space $\A$ is said {\it involutive}.

\vskip0.3cm

For $Q\in \A$ and for $A_1,A_2 \in \a$, let $R_Q(A_1,A_2)$ denote the equivalence class
in $\g/\a\oplus\b$ represented by $[A_1 + Q(A_1),A_2 + Q(A_2)]$. $R_Q$
defines a $\g/\a\oplus\b$-valued two-form on $\a$:
\[
 R_Q : (A_1,A_2)\in \a\times \a \mapsto R_Q(A_1,A_2)\in \g/\a\oplus\b.
  \]
If $R_Q$ vanishes identically, we may consider the
$\b$-valued two-form
$\tau_Q\in \b\otimes\Lambda^2(\a^\ast)$ defined by
\[
 \tau_Q(A_1,A_2) := [A_1 + Q(A_1),A_2 + Q(A_2)]_{\b} - Q\left([A_1 + Q(A_1),A_2 + Q(A_2)]_{\a}\right),
  \]
where, for $X\in \a\oplus\b$, $X_{\a}$ (resp. $X_{\b}$) denotes the $\a$ (resp. $\b$)
component of $X$.

\begin{defn}\label{def:tableau}
Let $\g$ be a Lie algebra and $\a, \b$ subspaces of $\g$ as above. By a
{\it tableau over} $\,\g$ we mean a subspace $\A$ of $\hom(\a,\b)$ 
(i.e., a tableau in the usual sense) such that:

\begin{enumerate}
\item $R_Q = 0$, for each $Q\in \A$;
\item $\tau_Q \in \text{Im}\, (\rho)$, for each $Q\in \A$;
\item $\A$ is involutive.
\end{enumerate}

\end{defn}

Tableaux of this type arise naturally in the case of semisimple Lie algebras.

\subsection{Tableaux over semisimple Lie algebras}\label{ss:Cartan-tableau}

Let $\g$ be a semisimple Lie algebra
endowed with its Killing form $\langle\,, \rangle$.
Let $\g=\g_0\oplus \g_1$ be a symmetric 
decomposition of $\g$, that is,
$\g_0$ is a subalgebra of $\g$ and $\g_1$ is a vector subspace such that
\[
 [\g_0,\g_0]\subset \g_0, \quad [\g_0,\g_1]\subset \g_1, \quad [\g_1,\g_1]\subset \g_0.
  \]
Assume that $\g/\g_0$ have rank $k$ and that
$\a$ be a maximal ($k$-dimensional) Abelian subspace of $\g_1$. Then the subspace $\g_1$ decomposes
as $\g_1 = \a\oplus \m$, where
\[
 \m = \a^\perp\cap \g_1 =\{X\in \g_1 \, : \, \langle X, A \rangle =0, \,\text{for all}\,\, A \,\in \a\}.
  \]
Further, let
\begin{gather}
 (\g_0)_\a = \{X\in \g_0 \,: \, [X,\a]=0\}, \\ (\g_0)_\a^\perp = \{X\in \g \,: \,
  \langle X, Y \rangle =0, \,\text{for all}\, \, Y \in (\g_0)_\a\},
   \end{gather}
and set $\b=\g_0\cap (\g_0)_\a^\perp$.
Then, for any regular\footnote{We recall that
$A\in \g_1$ is regular if the orbit at $A$ for the $\text{Ad}(G_0)$-action is a principal orbit.
If $A$ is regular, it is contained
in a maximal abelian subspace  $\a$ in $\g_1$ and $\text{Ad}(G_0)(\a)$ is open in $\g_1$ 
(cf. \cite{BDPT02, Ter97, TW03}).}
element $A\in \a$ with respect to the $\text{Ad}(G_0)$-action on $\g_1$, the maps
\[
 \text{ad}_A : \m \to \b, \quad \text{ad}_A : \b \to \m
  \]
are vector space isomorphisms. Moreover, the mapping
\[
 X\in  \m \mapsto -\text{ad}_X \in \hom(\a,\b)
  \]
is injective, and hence $\m$ can be identified with a linear subspace of $\hom(\a,\b)$.

\begin{prop}\label{Cartantableau}
If $\g$ is a semisimple Lie algebra and $\a$, $\b$, and $\m$ are defined as above,
then $\m$, regarded as a subspace of $\hom(\a,\b)$, is a tableau over $\g$.
It is referred to as a {\em Cartan tableau} over $\g$.
\end{prop}

\begin{proof}
First, we prove that $R_Q=0$, for every $Q= -\text{ad}_X$, $X\in\m$. For $A_1,A_2 \in \a$, we have
\[
 [A_1 + Q(A_1),A_2 + Q(A_2)] = \left[A_1 -[X,A_1],A_2 -[X,A_2]\right]= \left[[X,A_1],[X,A_2]\right].
  \]
By Jacobi's identity,
\[
 2\left[[X,A_1],[X,A_2]\right] = - \left[\left[[X,A_2],X\right],A_1\right] +
  \left[\left[[X,A_1],X\right],A_2\right].
   \]
Since $\left[[X,A_2],X\right]$ and $\left[[X,A_1],X\right]$ belong to $\left[[\g_1,\g_1],\g_1\right]\subset
[\g_0,\g_1]\subset \g_1$, we can write
\[
 \left[[X,A_2],X\right] = A' +B', \quad \left[[X,A_1],X\right] = A'' + B'',
  \]
for some $A',A'' \in \a$ and $B',B'' \in \m$. As a consequence,
\[
\left[\left[[X,A_2],X\right],A_1\right] = [B',A_1]\in \b, \quad
\left[\left[[X,A_1],X\right],A_2\right] = [B'',A_2]\in \b,
\]
which implies $[A_1 + Q(A_1),A_2 + Q(A_2)]\in \b$, for all $A_1,A_2\in \a$, and hence $R_Q=0$.

\vskip0.2cm

Next, we show that $\tau_Q\in \text{Im}\,(\rho)$, for every $Q=-\text{ad}_X$, $X\in\m$. From
the first part of the proof, we know that $[A_1 + Q(A_1),A_2 + Q(A_2)]\in\b$ and then
\[
 \tau_Q(A_1,A_2) =
  \frac{1}{2}\left(\left[\left[[X,A_1],X\right],A_2\right]- 
\left[\left[[X,A_2],X\right],A_1\right]\right),
   \]
that is,
\[
 \tau_Q(A_1,A_2) = \frac{1}{2}\left(\text{ad}_{M_Q(A_1)}(A_2) - \text{ad}_{M_Q(A_2)}(A_1) \right),
  \]
where $M_Q$ is the liner map defined by $M_Q : A\in \a \mapsto \left[[X,A],X\right]_{|\m}$
(here $\left[[X,A],X\right]_{\m}$ denotes the $\m$ component of $\left[[X,A],X\right]\in\g_1)$. This shows
that $\tau_Q\in \text{Im} (\rho)$.

\vskip0.2cm

Finally, we prove that $\m$ is involutive. Let $(A_1,\dots,A_k)$ be a basis consisting of regular
elements of $\a$. Then the flag
\[
 \ell : (0)\subset \text{span}\,\{A_1\} \subset \text{span}\,\{A_1,A_2\}\subset \cdots\subset
  \text{span}\,\{A_1,\dots,A_k\} =\a
   \]
is generic with respect to $\m$ and the characters are given by
\[
 s_1' = \text{dim}\, \m, \quad s'_j = 0 \quad (j=2,\dots,k).
  \]
For every $X\in \m$, we consider the linear map $\mu_X : \a \to \m$ defined by
\[
 [\mu_X(A_j), A_1] = [X,A_j], \quad \text{for all} \, A_j.
  \]
Then $\mu : X\in \m \mapsto \mu_X\in \hom(\a,\b)$ is a linear embedding such that $[\mu_X(A), B]=
[\mu_X(B), A]$, for every $X\in \m$ and every $A,B\in \b$. This shows that $\text{Im}\,(\mu) \subset
\m^{(1)}$. In particular, we get
\[
 \text{dim}\,\m \leq \text{dim}\,\m^{(1)} \leq s'_1 + 2s_2' + \cdots + ks_k' = \text{dim}\, \m.
  \]

\end{proof}

\section{The Pfaffian system associated with a tableau}\label{s:system-tableaux}

Let $\A\subset \hom(\a,\b)$ be a tableau over a finite dimensional Lie algebra $\g$ and let $G$ be
a connected Lie group with Lie algebra $\g$. We set $Y := G\times \A$ and refer to it as 
the configuration space.

A basis $(A_1,\dots,A_k,B_1,\dots,B_h, C_1,\dots,C_s)$ of $\g$ is said to be adapted to the tableau $\A$ if
\begin{equation}
 \a= \text{span}\,\{A_1,\dots,A_k\}, \quad \b = \text{span}\,\{B_1,\dots,B_h\}.
  \end{equation}
Moreover, we say that it is a {\it generic adapted basis} if the flag
\[
 (0) \subset \text{span}\,\{ A_1\} \subset\cdots \subset \text{span}\,\{A_1,\dots,A_k\} =\a
  \]
is generic with respect to $\A$.
For a generic adapted basis $(A_1$, $\dots,A_k$,$B_1$, $\dots,B_h$, $C_1$, $\dots$, $C_s)$, we let
\[
 (\alpha^1,\dots,\alpha^k,\beta^1,\dots,\beta^h,\gamma^1.\dots,\gamma^s)
  \]
denote its dual coframe. For a given basis
\[
 Q_\epsilon = Q^j_{\epsilon i} B_j \otimes \alpha^i \quad (\epsilon=1, \dots m)
  \]
of the tableau $\A$, we identify the
configuration space $Y$ with $G\times \R^m$ by
\begin{equation}
 (g,p^1,\dots,p^m) \in G\times \R^m \mapsto (g,  p^\epsilon Q_\epsilon) \in Y.
  \end{equation}

\begin{defn}
The exterior differential system associated with $\A$
is the Pfaffian system $(\mathcal{I}, \omega)$ on $Y$ such that $\I$ is generated,
as exterior differential ideal,  by
the linearly independent 1-forms
\begin{subequations}\label{system}
\begin{gather}
 \eta^j := \beta^j -
  p^\epsilon Q^j_{\epsilon i}\alpha^i \quad (j=1,\dots,h),\label{systema}\\
  \gamma^1, \dots, \gamma^s,  \label{systemb}
   \end{gather}
    \end{subequations}
with the independence condition
\begin{equation}\label{ind-cond}
\omega = \alpha^1 \wedge \cdots\wedge \alpha^k \neq 0.
\end{equation}
\end{defn}

The Pfaffian system $(\I,\omega)$ is given by a filtration of subbundles $I\subset J\subset T^\ast Y$
such that $J/I$ has rank $k$. The 1-forms $\eta^j$ $(j=1,\dots,h)$ and $\gamma^a$ $(a=1,\dots,s)$ span
the sections of $I$, while the $\alpha^i$ $(i=1,\dots,k)$ are sections of $J$ which project to a coframe of
$J/I$ at each point of $Y$.

\begin{thm}\label{thminvolutivity}
Let $\A$ be a tableau over a Lie algebra $\g$. Then, the exterior differential system
$(\I,\omega)$ on $Y$ associated with $\A$ is a linear Pfaffian system in involution.
In particular, the characters
$s'_j$ of $\A$ coincide with the Cartan characters $s_j$ of the system.
\end{thm}

\begin{proof}

Let  $V_k(\I,\omega)$ denote the set of all integral elements of
dimension $k$. We remind that $(\I,\omega)$ is in involution at
$y\in Y$ if there exists an ordinary integral element $(y,E) \in
V_k(\I,\omega)(y)$ and that $(\I,\omega)$ is in involution if it
is in involution at every point $y\in Y$. To prove the theorem, we
use Cartan's test of involution for linear Pfaffian systems.

The set of 1-forms $\{dp^\epsilon\} = \{dp^1,\dots,dp^m\}$
completes $\{\eta^j;\gamma^a;\alpha^i\}$ to a local coframe of $Y$
which is adapted to $I\subset J\subset T^\ast Y$. In terms of this
coframe, according to (1) of Definition \ref{def:tableau}, the
structure equations of $(\I,\omega)$ become
\begin{subequations}\label{str-system}
\begin{align}
 d\eta^j &\equiv -Q^j_{\epsilon\,i} dp^\epsilon\wedge\alpha^i + 
 \frac{1}{2}T^j_{i\,l}\alpha^i\wedge\alpha^l \mod \{I\}
 \quad &&(j=1,\dots,h),\label{str-systema}\\
  d\gamma^a &\equiv 0 \mod \{I\} \quad &&(a=1,\dots,s),  \label{str-systemb}
   \end{align}
    \end{subequations}
where
\[
 T^j_{i\,l}= - T^j_{l\,i}
  \]
and $\{I\}$ denote the algebraic ideal generated by $\{\eta^j;\gamma^a\}$.
From \eqref{str-system} it follows that the Pfaffian system $(\I,\omega)$ is linear, 
i.e., $dI \equiv 0 \mod \{J\}$.
This is equivalent to the
condition that the fibers of $V_k(\I,\omega)$ are affine
linear subspaces of $G_k(TY,\omega)$, with respect to the standard coordinates there. More specifically, an
integral element
at $y\in Y$ is defined by the equations
\begin{equation}\label{linear-eqs1}
 \eta^j=0, \quad \gamma^a =0, \quad dp^\epsilon = p^\epsilon_i\alpha^i,
  \end{equation}
for fiber coordinates $p^\epsilon_i$ which satisfy
\begin{equation}\label{linear-eqs2}
 -Q^j_{\epsilon\,i}(y)p^\epsilon_i + Q^j_{\epsilon\,i}(y)p^\epsilon_j - T^j_{i\,l}(y) =0.
   \end{equation}

Under a change of coframe
\[
 \tilde{\alpha}^i = \alpha^i, \quad \tilde{\eta}^j = \eta^j, \quad \tilde{\gamma}^a = \gamma^a,
  \quad \tilde{\pi}^\epsilon = dp^\epsilon - x^\epsilon_i\alpha^i,
   \]
the numbers $T^j_{i\,l}(y)$ transform to
\begin{equation}\label{torsion}
 \tilde{T}^J_{i\, l}(y) = -Q^j_{\epsilon\,i}(y)x^\epsilon_i + Q^j_{\epsilon\,i}(y)x^\epsilon_j - T^j_{i\,l}(y).
  \end{equation}
Then, using (2) of Definition \ref{def:tableau}, we may choose a
coframe of the form
\[
 \{\eta^j;\gamma^a;\alpha^i;\pi^\epsilon= dp^\epsilon -
   x^\epsilon_i\alpha^i\},
   \]
where $x^\epsilon_i$ are suitable smooth fuctions, so that the torsion term $T^j_{i\,l}$ in
\eqref{str-system} 
vanishes identically and the structure equations of
$(\I,\omega)$ take the form
\begin{align}
 d\eta^j &\equiv -Q^j_{\epsilon\,i} \pi^\epsilon \wedge\alpha^i 
 \mod \{I\} &&(j=1,\dots,h),\label{strrefined1}\\
  d\gamma^a &\equiv 0 \mod \{I\}
  &&(a=1,\dots,s).\label{strrefined2}
   \end{align}
Observe that (\ref{strrefined1}) and (\ref{strrefined2}) imply at
once that the integral elements
over $y\in Y$ form a linear space which is isomorphic to the first
prolongation $\A^{(1)}$ of the tableaux $\A$.

Next, set $\pi^j_i = -Q^j_{\epsilon\,i} \pi^\epsilon$ and
consider the tableau matrix 
\[
 \pi = \left(\begin{array}{c|c}
 
  \pi^j_i&{0}\\
\hline
{0}&{0}
\end{array}\right).
\]
If $A_1,\dots,A_k$
is generic with respect to $\A$, it is not difficult to show (cf. \cite{BCGGG}, p. 121) that
\begin{equation}\label{reduced-char}
 s'_1 + \cdots + s'_j = \left\{\begin{array}{l}
  \text{number of independent 1-forms}\\
   \text{$\pi^j_i$ in the first $j$ columns of $\pi$}
    \end{array} \right\},
     \end{equation}
where $s'_1,\dots,s'_k$ are the characters of the tableaux.
This shows that $s'_1,\dots,s'_k$ coincide with
the Cartan's characters $s_1,\dots,s_k$ of the differential system. Further,
since $\mathcal{A}$ is involutive and since
$\mathrm{dim}\,V_k(\I,\omega)(y)=\mathrm{dim}\,\A^{(1)}$, it
follows that
\[
 \mathrm{dim}\, V_k(\I,\omega)(y) = s_1+2s_2+\cdots +ks_k,\quad \text{for each} \quad y\in Y.
  \]
From Cartan's test of involution for linear Pfaffian system, we may then conclude that
 $(\mathcal{I},\omega)$ is involutive, with Cartan's characters $s'_1,\dots,s'_k$.
\end{proof}

\subsection{The $G/G_0$-exterior differential system}\label{ss:terng-sys}

Let $G/G_0$ be a semisimple symmetric space of rank $k$ and $\g = \g_0 \oplus \g_1$ be 
a Cartan
%symmetric
decomposition of $G/G_0$. 
%(this decomposition need not be a Cartan decomposition as the Killing
%form $\langle\,,\rangle$ may have signature on $\g_0$ and $\g_1$).
Let $\a$ be a maximal $k$-dimensional Abelian subspace in $\g_1$ and $(A_1,\dots,A_k)$
a basis for $\a$ consisting of regular elements with respect to the $\text{Ad}(G_0)$-action on $\g_1$.
According to Terng \cite{Ter97}, the {\it $G/G_0$-system} associated with the symmetric space $G/G_0$ is
the following system of PDEs for maps $V : U\subset \a \to \g_1\cap\a^\perp\,$:
\begin{equation}\label{ggzerosystem}
 [A_i, V_{x_j}]-[A_j, V_{x_i}]= [[A_i, V], [A_j, V]], \quad 1\leq i\not=j\leq k,
  \end{equation}
where $V_{x_i} = \frac{\partial V}{\partial x_i}$, being $x_i$ the coordinates
with respect to $(A_1,\dots,A_k)$.

\vskip0.2cm
Observe that $V : \a \to \g_1\cap\a^\perp$ is a solution of \eqref{ggzerosystem} if and only if
the $\g$-valued 1-form
\[
  \theta = \alpha + [\alpha,V] \in \Omega^1(\a)\otimes \g,
  \]
is flat, i.e., satisfies the Maurer--Cartan equation,
where $\alpha = \alpha^i\otimes A_i$ denotes the tautological 1-form on $\a$.

\vskip0.2cm
By a {\it generalized solution} of the $G/G_0$-system is meant a pair $(f,V) : M^k \to \a\times
(\g_1\cap\a^\perp)$ defined on a $k$-dimensional manifold $M^k$ such that:

\begin{itemize}

\item $f : M^k \to \a$ is a local diffeomorphism;

\item the $\g$-valued 1-form $df + [df,V] \in \Omega^1(\a)\otimes \g$ satisfies the Maurer--Cartan
equation.

\end{itemize}

\begin{remark}
Note that if $(f,V)$ is a generalized solution and $U\subset M^k$ is an open subset such
that $f_{|U} : U \to \a$ is a diffeomorphism onto its image, then
$V\circ (f_{|U})^{-1} : f(U)\subset \a \to \g_1\cap\a^\perp$
is a solution to the $G/G_0$-system.
\end{remark}

\begin{defn}
Let $(f,V)$ be a generalized solution of the $G/G_0$-system and let $\theta$ be the corresponding 1-form.
As $\theta$ satisfies the Maurer--Cartan equation, then there exists a map $g : M \to G$, uniquely
defined up left multiplication by an element of $G$, such that $\theta = g^{-1}dg$ (the pull-back by $g$ of 
the Maurer--Cartan form on $G$). $(g,V)$ is called
a {\it framed solution} of the $G/G_0$-system.
\end{defn}

For a map $(g,V) : M^k \to G\times \g_1\cap\a^\perp$, let $\alpha$ denote the $\a$ component of $g^{-1}dg$
with respect to the splitting $\g = \g_0 \oplus \a \oplus \g_1\cap\a^\perp$ of $\g$. Then $(g,V)$ is a
framed solution if and only if
\begin{itemize}
\item $g^{-1}dg=\alpha + [\alpha,V]$,

\item $\alpha^1\wedge\cdots \wedge \alpha^k \neq 0$.
\end{itemize}
If $M^k$ is simply connected, the corresponding generalized solution is given by $(f,V)$, where
$f : M^k \to \a$ is defined by $df =\alpha$ (note that $\alpha$ is a closed 1-form).

\vskip0.2cm

\begin{defn}
The {\it $G/G_0$-exterior differential
system} is defined to be the Pfaffian
system on the configuration space $Y=G\times \m$ associated with the Cartan tableau
$\m \subset \hom(\a,\b)$ .
\end{defn}

Next, we show that the framed solutions of the $G/G_0$-system can be described as integral manifolds of this
$G/G_0$-exterior differential system.

\vskip0.3cm 

Let $\n \subset \g_0$ be such that $\g_0 = \b \oplus \n$ and let
\[
 \theta  = \theta_\a + \theta_\m + \theta_\b + \theta_\n
  \]
be the decomposition of the Maurer--Cartan form of $G$ with respect to the splitting
\[
 \g= \a\oplus \m\oplus \b\oplus \n.
  \]
The Pfaffian system on $Y$ associated with the tableau $\m\subset \hom(\a,\b)$ is generated
by the Pfaffian equations
\[
 \theta_\b =[\theta_\a,V], \quad \theta_\m = \theta_\n =0,
  \]
with independence condition
\[
 \theta^1_\a\wedge\cdots \wedge \theta^k_\a\neq 0.
  \]
The integral manifolds of this system are smooth $k$-dimensional immersed submanifolds
\[
 (g,V) : M^k \to G\times (\g_1\cap\a^\perp)
  \]
which satisfy $g^{-1}dg= df + [df,V]$, where $f$ is a local diffeomorphism. This shows that the
integral manifolds coincide with the framed solutions of the $G/G_0$-system.

\begin{remark} 
According to the proof of Proposition \ref{Cartantableau} and
Theorem \ref{thminvolutivity}, it follows that the $G/G_0$-system
is in involution and its general solutions depend on $k$ functions
in one variable (cf. \cite{TW03}).
\end{remark}

\section{Surfaces in projective space and tableaux over $\sl(4,\R)$}\label{s:projective}

First, we briefly review some aspects of projective differential geometry of surfaces by
the method of moving frames (cf. \cite{AG, Bo, Ca, Fi}).

\subsection{The Wilczynski frame}

Let $G= \SL(4,\R)/\Z_2$ be the full group of projective transformations of the real projective
space $\P^3$. For each $g =(g_0, g_1, g_2, g_3) \in \SL(4,\R)$ we shall denote by $[g]$ its
equivalence class in $G$. The Lie algebra $\g$ of the projective group will be identified with
$\sl(4,\R)$ and the Maurer--Cartan form of $\g$ will be denoted by $\theta =(\theta^i_j) = g^{-1}dg$.

\vskip0.2cm

Let $f : M^2 \to \P^3$ be a smooth immersion of a connected surface. We shall still write
$\theta$ instead of $f^\ast(\theta)$ to denote the pull-back of the Maurer--Cartan form on $M$.
A {\it projective frame field}
along $f$ is a smooth map $g : U \to \SL(4,\R)$ defined on an open set $U$ of $M$ such that
$f(x) = [g_0(x)]$, for all $x\in U$. A projective frame along $f$ is of {\it first order} if
$\theta^3_0=0$. It easily seen that first order frames exist locally near any point of $M$.
Note that the 1-forms $\theta^1_0$ and $\theta^2_0$ define a coframe on $M$.
Differentiating $\theta^3_0=0$, it follows, by the structure equations and Cartan's Lemma, that there
exist smooth functions
$h_{ij} =h_{ji}$, $1\leq i,j \leq 2$, defined on $U$ and depending on the frame field, so that
\begin{equation}
 \theta^3_i = h_{ij}\theta^j_i, \quad i=1,2.
  \end{equation}
The vector valued quadratic form
\begin{equation}
 \phi_f = h_{ij}\theta^i_0 \theta^j_0 \otimes g_0
  \end{equation}
is independent of the choice of the first order frame and gives rise to a global cross section of
the vector bundle $S^2(M)\otimes K_f$, where $S^2(M)$ denotes the bundle of symmetric tensor of
type $(2,0)$ and $K_f$ is the canonical line bundle of the surface.

\begin{defn}
$\phi_f$ is known as the {\it Fubini quadratic form} of the surface. If $\phi_f$ has signature $(1,1)$,
then $f : M \to \P^3$ is said to be an {\it hyperbolic surface}.
\end{defn}

\begin{remark}
Most of the classical literature deals with surfaces of hyperbolic type. From now on
we shall restrict our consideration to such surfaces.
\end{remark}

Differentiating the equations $\theta^3_i = h_{ij}\theta^j_i$ ($i=1,2$) and applying the structure equations
and Cartan's Lemma, we have that
\begin{equation}
 -dh_{ij} + h_{ij}\theta^k_i + \frac{1}{2}h_{ij}(\theta^0_0 + \theta^3_3) = F_{ijk}\theta^k_0,
  \end{equation}
where the $F_{ijk}$ are smooth functions on $U$ symmetric in $i,j$ and $k$.
The symmetric cubic form defined by
\[
  F_{ijk}\theta^i_0\theta^j_0\theta^k_0 \otimes g_0
  \]
depends on the first order frame $g$ and transforms by
a section $\rho\cdot \phi_f$ of the line subbundle $S^1(M)\otimes \phi_f$ of $S^3(M)\otimes K_f$.

\begin{defn}
The equivalence class
\[
 \Psi_f = [F_{ijk}\theta^i_0\theta^j_0\theta^k_0 \otimes g_0]
  \]
in $S^3(M)\otimes K_f/S^1(M)\otimes \phi_f$ defines a global cross section, called
the {\it Fubini cubic form} of the surface. If $\Psi_f$ never vanishes, then 
$f : M^2 \to \P^3$
is said to be {\it generic of hyperbolic type}.
\end{defn}

\begin{remark}
The vanishing of the Fubini cubic form characterizes hy\-per\-bol\-ic quad\-rics in $\P^3$.
\end{remark}

By successive frame reductions one can prove the classical

\begin{thm}[Cartan, cf. \cite{Ca, Fi, AG}]\label{existswilz}
If $f : M^2 \to \P^3$ is a generic hyperbolic surface, then there
exists a unique projective frame field $[g] : M^2 \to G$ along $f$,
referred to as {\em the Wilczynski frame field},
such that
\begin{equation}\label{wil}
g^{-1}dg = (\theta^i_j)_{0\leq i,j\leq 3} = \left(\begin{array}{cccc}
  \theta^0_0 & \theta^0_1 & \theta^0_2 & \theta^0_3\\
    \theta^1_0 & \theta^1_1 & \theta^2_0 & \theta^0_2\\
 \theta^2_0 & \theta^1_0 & -\theta^1_1 & \theta^0_1\\
  0 & \theta^2_0 & \theta^1_0 & -\theta^0_0\\
\end{array}\right)
\end{equation}
and
\begin{equation}\label{indcondition}
{\theta^1_0\wedge \theta ^2_0} \neq 0.
\end{equation}
\end{thm}

\begin{remark}
In terms of the Wilczynski frame, the quadratic and cubic forms can be written, respectively, as:
\[
 \phi_f = \theta^1_0\cdot\theta^2_0 \otimes g_0, \quad \Psi_f = [(\theta^1_0)^3 + (\theta^2_0)^3 \otimes
   g_0].
   \]
\end{remark}

\vskip0.2cm
The non zero components of the Maurer--Cartan form $\theta$ of the Wilczynski frame are
$\theta^1_0$, $\theta^2_0$ (which define a coframe on $M^2$) and $\theta^0_0$, $\theta^1_1$, $\theta^0_1$,
$\theta^0_2$, $\theta^0_3$.
From the exterior differentiation of these forms and the structure equations, it follows that there
exist smooth functions $q_1,q_2$, $p_1,p_2$, and $r_1,r_2$, the {\it invariant functions of the surface},
such that
\begin{equation}\label{invfun}
 \left\{ \begin{array}{lll}
 \theta^0_0=-\frac{3}{2}q_1\theta^1_0 + \frac{3}{2}q_2\theta^2_0, &
  \theta^1_1= \frac{1}{2}q_1\theta^1_0+ \frac{1}{2}q_2\theta^2_0,\\
   \theta^0_1= r_1\theta^1_0 + p_2\theta^2_0,&
    \theta^0_2= p_1\theta^1_0 + r_2 \theta^2_0,\\
     \theta^0_3= r_2\theta^1_0 + r_1\theta^2_0,
      \end{array}\right.
       \end{equation}
In addition, the invariant functions must satisfy the following equations:
\begin{equation}\label{dteta1-2}
  d\theta^1_0=-q_2\theta^1_0\wedge\theta^2_0,\quad
   d\theta^2_0=-q_1\theta^1_0\wedge\theta^2_0,
    \end{equation}
\begin{equation}\label{dq1-dq2}
  \left\{ \begin{array}{rcl}
 dq_1\wedge \theta^1_0 - dq_2\wedge \theta^2_0 &=& \frac{2}{3}(p_1-p_2)\theta^1_0\wedge\theta^2_0,\\
   dq_1\wedge \theta^1_0 + dq_2\wedge \theta^2_0 &=& 2(1+q_1q_2-p_1-p_2)\theta^1_0\wedge \theta^2_0,
     \end{array}\right.
      \end{equation}
\begin{equation}\label{dr1-dr2}
  \left\{ \begin{array}{rcl}
   dr_1\wedge \theta^1_0 + dp_2\wedge \theta^2_0 &=& (2q_2r_1+3q_1p_2)\theta^1_0\wedge \theta^2_0,\\
    dp_1\wedge \theta^1_0+ dr_2\wedge \theta^2_0 &=& (2q_1r_2+3q_2p_1)\theta^1_0\wedge \theta^2_0,\\
     dr_2\wedge \theta^1_0+ dr_1\wedge \theta^2_0 &=& 4(q_1r_1+q_2r_2)\theta^1_0\wedge \theta^2_0.
       \end{array}\right.
        \end{equation}
Conversely, given a coframe $(\theta^1_0,\theta^2_0)$ on a simply connected surface $M$ and a set of
six real-valued
functions  $q_1,q_2$, $p_1,p_2$, $r_1,r_2$ satisfying \eqref{dteta1-2}, \eqref{dq1-dq2} and \eqref{dr1-dr2},
then there
exists a generic hyperbolic immersion $f : M \to \P^3$ with canonical coframe $(\theta^1_0,\theta^2_0)$
and invariant functions $q_1,q_2, p_1,p_2, r_1, r_2$.

\subsection{The Fubini--Cartan tableau and the corresponding exterior differential system}

Let $(\alpha^1,\alpha^2,\beta^1,\dots,\beta^5,\gamma^1.\dots,\gamma^8)$ be the basis of $\g^\ast$ defined
by
\begin{gather*}
 \alpha^1 = \theta^1_0, \quad \alpha^2 = \theta^2_0, \quad \beta^1=\theta^0_0, \quad \beta^2 =\theta^1_1,
  \quad \beta^3 = \theta^0_1, \quad \beta^4 = \theta^0_2, \quad \beta^5 = \theta^0_4,\\
   \gamma^1 = \theta^0_3, \quad \gamma^2 = \theta^2_0 -\theta^3_1, \quad \gamma^3 = \theta^1_0 -\theta^3_2,
   \quad \gamma^4 = \theta^1_1 + \theta^2_2, \\
   \gamma^5 = \theta^1_0 -\theta^2_1, \quad \gamma^6 = \theta^2_0 -\theta^1_2, \quad
   \gamma^7 = \theta^0_1 -\theta^2_3,
   \quad \gamma^8 = \theta^0_2 - \theta^1_3.
\end{gather*}
Let $(A_1,A_2,B_1,\dots,B_5, C_1,\dots,C_8)$ be its dual basis and set
\[
 \a = \text{span}\,\{A_1,A_2\}, \quad \b=\text{span}\,\{B_1,\dots,B_5\}.
  \]

\begin{defn}
The {\it Fubini--Cartan tableau} is the 6-dimensional subspace $\mathcal{W}\subset \hom(\a,\b)$
consisting of all elements $Q(q,p,r)$ of the form
\begin{align*}
 Q(q,p,r) &= q_1(-\frac{3}{2} B_1 + \frac{1}{2}B_2)\otimes \alpha^1
  +q_2(\frac{3}{2} B_1 + \frac{1}{2}B_2)\otimes \alpha^2 \\
&\qquad + p_1 B_4 \otimes \alpha^1 + p_2 B_3 \otimes \alpha^2
+ r_1(B_3 \otimes \alpha^1+ B_5\otimes \alpha^2) \\
 &\qquad   +  r_2(B_4 \otimes \alpha^2+ B_5\otimes \alpha^1),
\end{align*}
where $q=(q_1,q_2)$, $p=(p_1,p_2)$, $r=(r_1,r_2) \in \R^2$.

\end{defn}

%%%%%%%%%%%%%%%%%%%%%%%%%%%%%%%%%%%%%%%%%%

\begin{lemma}
$\mathcal{W}$ is a tableau over $\sl(4,\R)$.
\end{lemma}

\begin{proof} For any $Q \in \mathcal{W}$, we compute
\begin{multline*}
 [A_1 + Q(A_1),A_2 + Q(A_2)] = q_2A_1 + q_1A_2 +(p_1-p_2)B_1 + (p_2+p_1-1)B_2 \\- (q_2r_1 +2q_1p_2)B_3
  - (q_1r_2 +2q_2p_1)B_4 -3(q_1r_1 + q_2r_2)B_5.
  \end{multline*}
This implies that $[A_1 + Q(A_1),A_2 + Q(A_2)]\in \a\oplus \b$, for each $Q\in \mathcal{W}$, $A,B\in \a$.
A direct computation shows that the first prolongation $\mathcal{W}^{(1)}$ of $\mathcal{W}$ is
spanned by
\begin{gather*}
 Q_1 \otimes \alpha^1, \quad Q_2 \otimes \alpha^2,
  \quad Q_3 \otimes \alpha^1,
   \quad Q_4 \otimes \alpha^2,\\
    Q_4 \otimes \alpha^1 + Q_5 \otimes \alpha^2,
    \quad Q_5 \otimes \alpha^1 + Q_6 \otimes \alpha^2, \quad Q_6 \otimes \alpha^1 + Q_3 \otimes \alpha^2.
     \end{gather*}
Thus $\text{dim}\,\mathcal{W}^{(1)} = 7$ and the map
$\rho : \hom(\a,\mathcal{W}) \to \b\otimes \Lambda^2(\a^\ast)$ is surjective\footnote{in fact,
$\text{dim}\,\hom(\a,\b) =12$,
$\text{dim}\, (\b\otimes \Lambda^2(\a^\ast))=5$ and $\text{dim}\,\text{Ker}\,\rho =
\text{dim}\,\mathcal{W}^{(1)} = 7$.}.

Next, consider a generic flag $\ell \in \mathcal{F}(\a)$ of the form
\[
 \ell : (0) \subset \text{span}\,\{\lambda A_1 + \mu A_2\} \subset \a,
  \]
for $\lambda$, $\mu\in \R$, $\lambda\mu \neq 0$. The linear space consisting of all $Q\in \mathcal{W}$
such that $Q_{|\text{span}\{\lambda A_1 + \mu A_2\}}=0$ is characterized by the equations
\[
  q_1=q_2=0,\quad \lambda p_1 + \mu r_2 = \mu p_2 + \lambda r_1 = \mu r_1 + \lambda r_2 =0.
  \]
This implies that $\text{dim}\, \mathcal{W}_1 =1$ and therefore
\[
 s'_1 = \text{dim}\, \mathcal{W} -1= 5, \quad s'_2 =6-s'_1 = 1.
  \]
In conclusion, $\text{dim}\,\mathcal{W}^{(1)} = 7 = s_1' +2s'_2$, that is, $\mathcal{W}$ is involutive.
\end{proof}

%%%%%%%%%%%%%%%%%%%%%%%%%%%%%%%%%%%%%%%%%%%%%%%

\begin{remark}
The exterior differential system defined by the Fubini--Cartan tableau is the Pfaffian differential ideal
generated by the 1-forms
\[
 \gamma^1.\dots,\gamma^8, \eta^1, \cdots, \eta^5,
  \]
where $\eta^1, \cdots, \eta^5$ are given by
\begin{equation}\label{eta}
 \left\{ \begin{array}{lll}
 \eta^1 = \beta^1 + \frac{3}{2}q_1\theta^1_0 - \frac{3}{2}q_2\theta^2_0, &
  \eta^2 = \beta^2 - \frac{1}{2}q_1\theta^1_0 - \frac{1}{2}q_2\theta^2_0,\\
   \eta^3 = \beta^3 - r_1\theta^1_0 - p_2\theta^2_0,&
    \eta^4 = \beta^4 - p_1\theta^1_0 - r_2 \theta^2_0,\\
     \eta^5 = \beta^5 - r_2\theta^1_0 - r_1\theta^2_0,
      \end{array}\right.
       \end{equation}
with independent condition $\theta^1_0\wedge \theta^2_0 \neq 0$.
From this we infer that the integral manifolds of the system are the 2-dimensional
submanifolds
\[
 (g;q,p,r) : M^2 \to G\times \mathcal{W} \cong G\times \R^6
  \]
such that
\begin{itemize}
\item $f=[g_0] : M^2 \to \P^3$ is a generic hyperbolic surface;

\item $g : M^2 \to G$ is the Wilczynski frame along $f$;

\item $q_1,q_2,p_1,p_2, r_1,r_2 : M^2 \to \R$ are the invariant functions of $f$.
\end{itemize}

\end{remark}

\subsection{The projective Gauss map}

 We recall that the {\it Lie quadric} $Q(x)$ of the surface $f : M \to \P^3$
at $f(x) = [y^0:y^1:y^2:y^3]$ is defined, with respect to the Wilczynski frame $[g(x)]$, by the equation
\[
 y^0y^3 - y^1y^2 =0.
  \]
Thus we may represent $Q(x)$ by means of the quadratic form
\begin{equation}
\chi_f(x) = g^0(x)g^3(x) - g^1(x)g^2(x), \quad \text{for each} \quad x\in M.
 \end{equation}
Let $Q_{(2,2)}$ denote the pseudo-Riemannian symmetric space
\[
 Q_{(2,2)} = \SL(4,\R)/\SO(2,2).
  \]
Then the 2-parameter family of Lie quadrics can be viewed as an immersion
\[
\chi_f : M\to Q_{(2,2)}, x \mapsto \chi_f(x).
\]
This map is referred to as the {\it projective Gauss map} of the immersion $f$.

\vskip0.2cm
The Lie quadric $\chi_f(x)$ is parametrized by the map
\[
 \hat{f}_x : (\lambda,\mu) \mapsto g_0(x) + \lambda g_1(x) + \mu g_2(x) + \lambda \mu g_3(x).
  \]
Thus the projective tangent space of $\chi_f(x)$ at the point $\hat{f}_x(\lambda,\mu)$ is the projective
plane spanned by the following three linearly independent vectors
\begin{equation}
 \left\{\begin{array}{ll}
  \hat{f}_x(\lambda,\mu) = g_0(x) + \lambda g_1(x) + \mu g_2(x) + \lambda \mu g_3(x),\\
   \partial_\lambda\hat{f}_x(\lambda,\mu) = g_1(x) + \mu g_3(x),\\
     \partial_\mu\hat{f}_x(\lambda,\mu) = g_1(x) + \lambda g_3(x).
      \end{array}\right.
       \end{equation}
Accordingly, the envelopes of the projective Gauss map are maps $h : M \to \R^4$ of the form
\[
 h =  g_0 + \lambda g_1 + \mu g_2 + \lambda \mu g_3,
  \]
where $\lambda, \mu$ are smooth functions, determined by the tangency condition
\[
 dh\wedge\hat{f}_x \wedge \partial_\lambda\hat{f}_x \wedge \partial_\mu\hat{f}_x =0.
  \]
Using the Maurer--Cartan equations, we compute
\[
 dh\wedge\hat{f}_x \wedge \partial_\lambda\hat{f}_x \wedge \partial_\mu\hat{f}_x=
  (\lambda^2\theta^1_0 + \mu^2\theta^2_0 -\lambda^2\mu^2\theta^0_3)\,
   g_0\wedge g_1\wedge \wedge g_2 \wedge g_3.
    \]
As $\theta^0_3 =r_2\theta^1_0 +r_1 \theta^2_0$, then $\lambda$ and $\mu$ must be solutions of the
system
\[
  \lambda^2(1-\mu^2r_2) =0, \quad \mu^2(1-\lambda^2r_1)=0.
  \]
If we allow the envelopes to be imaginary, we have the classical result

\begin{prop}[cf. \cite{Bl2,Bo,Fi}]
The projective Gauss map of a nondegenerate hyperbolic surface $f : M \to \P^3$ has in
general five envelopes; out of these, one is the surface itself, while the other
four are given by
\begin{equation}\label{envelopes}
  h_{(\epsilon,\eta)} = g_0 + \epsilon \sqrt{r_2} g_1 + \eta\sqrt{r_1} g_2 +
   \epsilon\eta\sqrt{r_1r_2} \mu g_3,
     \end{equation}
where $\epsilon,\eta =\pm 1$. There are two exceptional cases:

\begin{itemize}

\item if one of the two invariant functions $r_1$ or $r_2$ vanishes identically, say $r_2=0$,
$r_1\neq 0$, then the projective Gauss map has three distinct envelopes: the surface itself and
$h_\pm = g_0 \pm \sqrt{r_1}\,g_2$;

\item if $r_1 = r_2 =0$, then there are only two envelopes: the surface itself and its {\em dual surface}
$\check{f} =[g_3]$.

\end{itemize}

\end{prop}

\begin{remark}
We have real envelopes if and only if $r_1>0$ and $r_2 >0$.
\end{remark}

\subsection{Godeaux--Rozet surfaces and the corresponding tableau}

{\it Godeaux--Rozet surfaces} are characterized by the fact that either $r_1=0$ or else $r_2=0$.
So the projective Gauss map of a Godeaux--Rozet surface has less than five distinct envelopes.
They can be viewed as integral submanifolds of
the Fubini--Cartan differential system restricted to one of the two submanifolds
\[
 Y'_{GR} = \{(g,q,p,r) \in Y \, : \, r_2=0 \}, \quad Y''_{GR} = \{(g,q,p,r) \in Y \, : \, r_1=0 \}
  \]
of $Y$. We consider $Y'_{GR}$, the other case being analogous. It easily seen that $Y'_{GR}$ is
the configuration space of the 5-dimensional tableau over $\sl(4,\R)$
\[
 \mathcal{W}_{GR} = \left\{Q(q,p,r) \in \mathcal{W} \, : \, r_2=0 \right\},
  \]
referred to as the {\it Godeaux--Rozet tableau}. Notice that the restriction to
$Y'_{GR}$ of the Fubini--Cartan differential system is exactly the system associated
with $\mathcal{W}_{GR}$. Its integral manifolds are canonical lifts of Godeaux--Rozet surfaces
in 3-dimensional projective space. It is easily checked that $\mathcal{W}_{GR}$ is
involutive and that its characters are
\[
  s'_0=13, \quad s'_1 =5, \quad s'_2 = 0.
   \]
Thus a general Godeaux--Rozet surface depends on five arbitrary functions in one variable.

\subsection{Demoulin surfaces and the corresponding tableau}

{\it Demoulin surfaces} are characterized by the fact that both $r_1=0$ and $r_2=0$.
So the projective Gauss map of a Demoulin surface has only two distinct envelopes.
Demoulin surfaces can be viewed as integral submanifolds of
the Fubini--Cartan differential system restricted to the submanifold
\[
 Y_{D} = \{(g,q,p,r) \in Y \, : \, r_1 = r_2=0 \}.
  \]
$Y_{D}$ is the configuration space of the tableau over $\sl(4,\R)$
\[
 \mathcal{W}_{D} = \left\{Q(q,p,r) \in \mathcal{W} \, : \, r_1=r_2=0 \right\},
  \]
referred to as the {\it Demoulin tableau}. The restriction to
$Y_{D}$ of the Fubini--Cartan differential system is exactly the system associated
with $\mathcal{W}_{D}$, whose integral manifolds are canonical lifts of Demoulin surfaces
in 3-dimensional projective space. It is an easy exercise to check that $\mathcal{W}_{D}$ is
involutive and that its characters are
\[
  s'_0=13, \quad s'_1 =4, \quad s'_2 = 0.
   \]
Thus a general Demoulin surface depends on four arbitrary functions in one variable.

\subsection{Asymptotically-isothermic surfaces}

A surface $f : M\to \P^3$ is called {\it asymptotically-isothermic} if the 3-web
defined by the families of curves satisfying the Pfaffian equations
\[
 \theta^1_0 = \theta^2_0 = \theta^1_0 + \theta^2_0 =0
  \]
is flat. Since the connection of this web is given by $\zeta_w=q_1\theta^1_0 -q_2\theta^2_0$,
we deduce that asymptotically-isothermic surfaces are characterized by the equation $p_1 -p_2 =0$.
The cor\-re\-spond\-ing 5-di\-men\-sional tableau is given by
\[
 \mathcal{W}_{AI} = \left\{Q(q,p_1, p_1,r)\in \mathcal{W}  \, : \, t \in \R, \, q,r \in\R^2 \right\}.
  \]
This tableau is involutive and its characters are 
\[
  s'_0=13, \quad s'_1 =5, \quad s'_2 = 0.
   \]
Therefore, asymptotically-isothermic surfaces depend on five arbitrary functions in one variable.

%%%%%%%%%%%%%%%%
%\begin{defn}
%A surface $f : M\to \P^3$ is said of {\it constant curvature} $c$ if the Fubini 
%quadratic form
%$\theta^1_0\cdot \theta^2_0$ has constant sectional curvature $c$. 
%The connection form of the
%Fubini metric is given by $\zeta_m=-q_1\theta^1_0 -q_2\theta^2_0$. 
%Thus $f$ has constant curvature
%$c$ if and only if its invariant functions $p_1$ and $p_2$ satisfy the constraint
%$p_1 +p_2 = 1 +\frac{c}{2}$
%\end{defn}
%%%%%%%%%%%%%%%

\begin{remark}

More generally, one could consider the class of surfaces whose invariant functions
$p_1$ and $p_2$ satisfy a linear relation, that is, are expressed by 
\begin{equation}\label{p1tp2t}
  p_1(t) = t \cos{a} + b_1, \quad p_2(t) = t \sin{a} + b_2,
   \end{equation}
where $a$, $b_1$, $b_2$ are real constants. 
This class of surfaces may be viewed as an analog of linear Weingarten surfaces and include, 
as special examples, the surfaces of constant curvature with respect to Fubini's quadratic form
$\theta^1_0\cdot \theta^2_0$. In fact, a surface $f : M\to \P^3$ has constant curvature
$c$ if and only if $p_1 +p_2 = 1 +\frac{c}{2}$.
Next, con\-sid\-er the 5-di\-men\-sional {\it affine}
tableau
\[
 \mathcal{W}_{(a,b_1,b_2)} = \left\{Q(q,p_1(t), p_2(t),r)  \, : \, t \in \R, \, q,r \in\R^2 \right\}.
  \]
Notice that $\mathcal{W}_{(a,b_1,b_2)}$ is
an affine subspace and as such, properly speaking, it does not fit into our scheme. 
Nevertheless, it is possible to develop a similar theory also for affine tableaux
over Lie algebras, and then surfaces satisfying \eqref{p1tp2t} can be viewed as integral 
submanifolds of the Pfaffian system on $G\times \mathcal{W}_{(a,b_1,b_2)}$ associated with 
$\mathcal{W}_{(a,b_1,b_2)}$.
Also in this case, a direct computation yields that $\mathcal{W}_{(a,b_1,b_2)}$ is involutive
and that its characters are
\[
  s'_0=13, \quad s'_1 =5, \quad s'_2 = 0.
   \]
Thus the surfaces whose invariant functions $p_1$ and $p_2$ satisfy a linear condition
(including  asymptotically-isothermic and constant curvature surfaces) depend on five arbitrary 
functions in one variable.

\end{remark}

\bibliographystyle{amsalpha}

\end{document}